# Le problème intensionnel de la quantification chez Hilbert

**Jean PETITOT**

EHESS et CREA (Ecole Polytechnique), Paris

petitot@ehess.fr

2002

**Abstract.**

The paper (in French) presents a survey of Hilbert's Epsilon operator focusing on the intensional aspects of its semantics. It comments on some epistemological problems, from Albert Lautman in the 1930s to John Bell, Grigori Mints, Barry Hartley Slater, Richard Zach or Edward Zalta in the 1990s.

## I. INTRODUCTION

### 1. L'importance du formalisme hilbertien

Cet article est une mise à jour d'anciennes réflexions des années 1970-1980 sur le formalisme introduit par Hilbert, Ackermann et Bernays pour maîtriser la quantification sur des domaines infinis. Il me semble en effet que ce formalisme – dit ε-calculus – est d'une importance philosophique remarquable. Je l'ai utilisé pour la première fois en 1975 lors de la rencontre sur l'apprentissage organisée entre Chomsky et Piaget au *Centre de Royaumont* par Massimo Piatelli. Je l'ai ensuite utilisé dans mon article *Infinitesimale* de 1979 dans l'*Enciclopedia Einaudi* pour clarifier en termes de logique philosophique le statut apparemment paradoxal des infinitésimales leibniziennes tel qu'il est formalisé par l'analyse non standard.[1] Puis j'ai donné plusieurs cours sur ce thème à l'Université de Bologne chez Umberto Eco en 1980 en présence de Marco Santambrogio qui a développé plus tard ces idées dans ses propres travaux et en particulier dans son ouvrage de 1992 *Forma e oggetto*. D'après le témoignage d'Alberto Peruzzi, ses discussions avec Santambrogio sur le cours de Bologne ont joué un rôle important dans ses propres réflexions sur la logique des descriptions définies.² Ces réflexions le conduisirent à entrer alors en contact avec John Bell et à poser le problème d'une formalisation de l'opérateur de Hilbert dans le cadre de la logique intuitionniste interne à un topos.

---

[1] Cf. aussi Petitot [1999].

**2.     Le problème de la quantification transfinie**

Le problème fondamental de la quantification qui se trouve à l'origine des travaux de Hilbert et de son école est le suivant. Une quantification existentielle $\exists x\, F(x)$ équivaut intuitivement à une disjonction $F(a_1) \vee \ldots \vee F(a_n) \vee \ldots$ où les $a_i$ sont les éléments du domaine de définition $dom(F)$ de $F$. De même, une quantification universelle $\forall x\, F(x)$ équivaut intuitivement à une conjonction $F(a_1) \wedge \ldots \wedge F(a_n) \wedge \ldots$ On pourrait donc croire que l'on peut éliminer les formules quantifiées au profit de telles disjonctions et conjonctions. Mais si le domaine $dom(F)$ de $F$ est *infini* alors $\vee$ et $\wedge$ deviennent des opérateurs *transfinis* et la vérification des énoncés devient impossible.

Le problème est incontournable. En effet d'après le théorème de Löwenheim-Skolem, toute théorie du premier ordre de domaine infini admet des modèles non standard.³ Prenons par exemple le cas d'une théorie s'exprimant dans le langage de l'arithmétique. Dans un modèle non standard $\mathbb{N}^*$ de $\mathbb{N}$ il existe des entiers $\nu$ infinis. Supposons alors que l'on sache démontrer pour tout $n \in \mathbb{N}$ l'énoncé $F(n)$. Cela n'impliquera pas en général l'énoncé universel $\forall n\, F(n)$. En effet cet énoncé étant du premier ordre est aussi valable dans $\mathbb{N}^*$, or il n'existe aucune raison pour que $F(\nu)$ soit vrai d'entiers infinis $\nu$. $F$ peut n'être valable que pour des entiers "vraiment" finis $n$. Or l'ensemble $\mathbb{N}$ des entiers "vraiment" finis *n'est pas* caractérisable au premier ordre comme sous-ensemble de $\mathbb{N}^*$, c'est un sous-ensemble "externe" et non pas "interne" de $\mathbb{N}^*$. Cela est très intuitif : $\mathbb{N}$ et $\mathbb{N}^*$ ayant la même théorie par définition des modèles non standard, si $\mathbb{N}$ était un sous-ensemble interne de $\mathbb{N}^*$, il y aurait un sous ensemble $\mathbb{N}°$ de $\mathbb{N}$ qui serait à $\mathbb{N}$ ce que $\mathbb{N}$ est à $\mathbb{N}^*$ et représenterait des entiers en quelque sorte "hyperfinis". Mais il est évident que, quelle que soit la définition envisagée, si $n$ est "hyperfini", $n+1$ le sera aussi et donc $\mathbb{N}° = \mathbb{N}$. Or $\mathbb{N} \neq \mathbb{N}^*$.

**3.     Les descriptions définies de Russell**

Pour approfondir le problème d'une possible élimination des quantificateurs à portée infinie, Hilbert s'est inspiré des descriptions définies introduites par Russell. Dans le cas où $F(x)$ satisfait les propriétés d'existence et d'unicité des $x$ tels que $F(x)$, Russell introduisit le symbole $\iota_x F(x)$ pour symboliser "*le x* tel que $F(x)$", ce qui conduit à définir dans ce cas la quantification existentielle par la formule : $\exists x\, F(x) \Leftrightarrow F(\iota_x F(x))$.

---

² Cf. Peruzzi [1983] et [1989].

³ Pour des précisions, cf. Petitot [1979] et sa bibliographie.

## II. DEFINITIONS GENERALES DE L'OPERATEUR DE HILBERT

### 1. ε-termes et quantification existentielle

Soit $F(x)$ un prédicat (unaire) quelconque. De façon purement syntaxique, Hilbert introduit un symbole *intensionnel d'individu* — dit *ε–terme* — noté $\varepsilon_x F(x)$ ou $\varepsilon_F$ par commodité, où $x$ devient une variable liée. $\varepsilon_x F(x)$ possédant le *même* type logique que celui des arguments $x$ de $F$, l'expression $F(\varepsilon_F)$, bien qu'auto-référentielle, est pourtant bien formée. Plus généralement, pour toute formule ouverte $G(y)$ à une variable libre $y$ de même type que $x$ (en particulier pour tout autre prédicat unaire), $G(\varepsilon_F)$ est une formule fermée, c'est-à-dire un énoncé. $\varepsilon_x F(x)$ est un symbole complet (fermé, saturé). Mais l'opérateur ε permet de construire des symboles incomplets. Si $F(x,y)$ est un prédicat binaire, alors $\varepsilon_x F(x,y)$ est un prédicat unaire de $y$ et l'on peut former l'ε-terme $\varepsilon_y \varepsilon_x F(x,y)$, etc.

L'existence et l'identité des *ε*-termes ainsi définis sont purement *syntaxiques*. Il est essentiel d'insister sur le fait que $\varepsilon_F$ est *un symbole d'individu et non pas de variable*. Sémantiquement (au sens banal de "sémantique"), $\varepsilon_F$ représente — et en quelque sorte "hypostasie" — l'*idée* d'un individu satisfaisant $F$, et cela que l'extension $X_F$ de $F$ soit non vide ou non. C'est un individu *idéal* représentant l'idée d'un individu qui satisfait $F$, autrement dit une *idée in individuo*. Comme le dit John Bell,

> "We may think of $\varepsilon_F$ as meaning an ideal object associated with $F$: all one knows about it is that, if anything satisfies $A$, it does."[4]

Ou encore, comme l'expliquent Urs Egli et Klaus von Heusinger (1995, p. 13), un ε-terme $\varepsilon_F$

> "denotes the most salient individual that has the property $F$"

où la saillance est définie par une hiérarchie de saillance au sens de Lewis (1979).

Hilbert *définit* alors la *quantification existentielle* par l'équivalence syntaxique :

(∃)    $\exists x\, F(x) \equiv F(\varepsilon_F)$.

---

[4] Bell [1993], p. 1.

Informellement interprétée, cette équivalence signifie que le fait qu'il existe un individu satisfaisant $F$ équivaut au fait que l'idée *in individuo* d'un individu satisfaisant $F$ satisfait effectivement $F$. Il s'agit donc d'un critère de *consistance de signification*. Par l'adjonction d'individus idéaux à l'univers considéré, il devient par conséquent possible de ramener des formules avec quantificateurs à des formules de forme strictement plus simple, *sans* quantificateurs. Il y a là une opération qui, au-delà de son intérêt proprement logique, est, comme nous l'affirmons depuis longtemps,[5] d'une grande portée philosophique.

Appelons *subsistance*, pour ne pas la confondre avec l'existence formalisée par le quantificateur $\exists$, l'existence purement symbolique et syntaxique des $\varepsilon$-termes.[6] L'équivalence ($\exists$) signifie que *l'existence équivaut à un redoublement en consistance (en cohérence sémantique) d'une subsistance*. L'existence est l'*effet* du fait que la subsistance symbolique $\varepsilon_F$ puisse posséder une référence consistante conforme à son contenu idéal. Elle est donc un type particulier d'*auto-référence*.

Au sens que possède l'opposition entre syntaxe et sémantique en théorie logique des modèles (où la sémantique s'identifie à la dénotation), qu'en est-il de l'interprétation des $\varepsilon$-termes ? A supposer que l'on ait interprété $F(x)$ dans un modèle $\mathfrak{M}$ du langage $\mathfrak{L}$ dont $F(x)$ est une formule, comment doit-on interpréter $\varepsilon_F$ ? Comme nous le verrons, de nombreuses discussions ont eu lieu à ce propos. Elles tournaient toutes, sans le reconnaître au départ, autour du fait que les $\varepsilon$-termes $\varepsilon_F$ sont des entités *intensionnelles*.

Comme un *déictique* d'une langue naturelle, $\varepsilon_F$ est en quelque sorte un *symbole-index* — i.e. une entité intensionnelle de nature *pragmatique* et *dépendante du contexte* — dont l'identité syntaxique est parfaitement bien définie (aspect symbole) mais dont l'identité sémantique (la dénotation) est au contraire indéfinie (aspect index). C'est pourquoi, avant que les logiques intensionnelles et la pragmatique formelle n'aient investigué ce genre d'entités, on interprétait sémantiquement les $\varepsilon$-termes comme des *opérateurs de choix*. Soit $X_F$ l'extension de $F$ et supposons $X_F \neq \emptyset$, i.e. supposons que $\exists x\, F(x)$ soit valide. Sémantiquement, $\varepsilon_F$ fait plus que simplement dénoter un certain élément $a \in X_F$ (ce qui serait le cas si $\varepsilon_F$ était simplement un symbole de constante). Il *sélectionne* un élément $a \in X_F$. La définition ($\exists$) devient alors la reformulation de la règle d'introduction du quantificateur existentiel (règle GE de généralisation existentielle) du calcul des prédicats :

---

[5] Nos premiers travaux sur cette question remontent au milieu des années 70 et concernaient l'utilisation du formalisme hilbertien pour invalider un soit disant "théorème" de Jerry Fodor ayant servi de base à la critique chomskyenne des théories associationnistes de l'apprentissage et à la justification de thèses innéistes.

[6] Nous empruntons cette dénomination à la théorie des états de choses (Sachverhalten, states of affairs). Cf. Petitot [1985].

(GE)    $F(a) \Rightarrow \exists x\, F(x)$.

En tant qu'*opérateur de choix* relatif à l'extension $X_F$, $\varepsilon_F$ peut être considéré comme un *type générique* dont les spécialisations sont les éléments $a \in X_F$. Ce simple fait suffit à en laisser soupçonner la pertinence cognitive et sémio-linguistique.

## 2. Le contexte philosophique : ontologie formelle et objets abstraits

Les constructions logiques hilbertiennes sont en profonde résonnance avec certains des plus anciens et des plus fondamentaux problèmes de la logique philosophique depuis le Moyen-Âge, ceux concernant l'ontologie formelle des *objets abstraits* idéaux et génériques ainsi que l'existence comme prédicat. L'histoire métaphysique en est immense puisque c'est celle du conflit réalisme / nominalisme en matière d'universaux.

Des nominalistes médiévaux à Berkeley et de Berkeley à Husserl, la plupart des philosophes ont refusé d'accepter des individus abstraits (i.e. des individus pouvant individuer une idée générale) pour la raison qu'en tant qu'entités abstraites se sont des objets incomplets alors qu'en tant qu'individus se sont nécessairement des objets complets d'après le tiers exclu (pour toute propriété $P$ ils doivent satisfaire $P$ ou $\neg P$). L'exemple classique est celui de l'absurdité du "triangle général" qui devrait être à la fois rectangle et équilatéral. C'est pour répondre à cette question que Kant inventa le concept de *schème* comme objet générique (en fait comme contenu noématique au sens de Husserl) associé à un concept et donnant la règle de construction des référents de ce concept.

Au XXème siècle cette question centrale a été relancée entre autres par Alexius Meinong et a fait l'objet d'un débat ayant impliqué tous les philosophes et logiciens intéressés par l'ontologie formelle, de Husserl à Quine en passant par l'école polonaise (Lesniewski, etc.). On trouvera un excellent panorama du débat dans l'ouvrage de Frédéric Nef (1998) sur *L'objet quelconque*.[7]

Comme nous le verrons plus bas, la réponse à ces difficultés se trouve bien du côté du statut pragmatique des ε-termes comme symboles-index : les ε-termes sont des individus en quelque sorte "intensionnellement décomplétés". Ce ne sont pas des entités abstraites mais des individus qu'une décomplétion intensionnelle (indexicale et modale) permet de considérer comme des intermédiaires entre des symboles de variables et des objets concrets.

---

[7] Pour une introduction à l'ontologie formelle on pourra aussi consulter Poli [1992].

Ce point a été particulièrement bien approfondi par le spécialiste d'ontologie formelle et de logique intensionnelle Edward Zalta qui défend l'introduction d'objets idéaux pour tous les types logiques et a consacré de nombreux travaux au prédicat d'existence, qu'il note E!$a$,[8] et au prédicat d'abstraction, qu'il note A!$a$ ("$a$ est un objet abstrait"). Son idée est que les individus génériques idéaux "*encodent*" des propriétés au lieu de les exemplifier. Ce ne sont pas des "supports" de propriétés mais des encodeurs de règles déterminantes et ils conduisent à définir *deux* types différents de *prédication*. Ce sont des "concepts-objets" (ce que nous appelons des idées *in individuo*) qui sont "constitués" par les descriptions qui les déterminent sans pour autant y satisfaire nécessairement. Ils ne sont donc complets que par rapport à leurs propriétés déterminantes et restent incomplets par rapport à leurs autres propriétés possibles. Cet écart entre "constitution" syntaxique et "satisfaction" sémantique est d'ailleurs une façon de formuler leur statut de symbole-index. Dans un article comparant la théorie de ces objets chez Ernst Mally (le disciple de Meinong lui ayant succédé à sa chaire de Gratz) et Husserl, E. Zalta (1998) conclut que, phénoménologiquement parlant, ce sont des contenus noématiques intentionnels

## 3.   ε-termes zéro et quantification universelle

Qu'en est-il maintenant d'un jugement non plus existentiel mais *universel*, du type $\forall x\, F(x)$ ? En logique *classique* (mais *pas* en logique *intuitionniste*, nous y reviendrons à la fin de cette étude), on a les équivalences suivantes :

$\forall x\, F(x) \Leftrightarrow \neg \exists x\, \neg F(x)$ (dualité $\forall \leftrightarrow \exists$ à travers la négation);[9]
$\exists x\, \neg F(x) \Leftrightarrow \neg F(\varepsilon_{\neg F})$ (d'après ($\exists$) appliquée à $\neg F$);
$\neg \exists x\, \neg F(x) \Leftrightarrow \neg \neg F(\varepsilon_{\neg F})$;
$\neg \neg F = F$ (double négation classique);

et donc, en définitive :

($\forall$)    $\forall x\, F(x) \Leftrightarrow F(\varepsilon_{\neg F})$.

---

[8] Le philosophe de la logique italien Giulio Preti a aussi utilisé le formalisme hilbertien pour définir un "prédicat d'existence" E! pour les ε-termes. Dans un inédit du 16 octobre 1955 intitulé *Ricerche ontologiche* il discute les équivalences E! $\varepsilon_x F(x) \Leftrightarrow F(\varepsilon_x F(x)) \Leftrightarrow \exists x\, F(x)$.

[9] $\neg F$ est la négation de $F$.

Il faut insister sur le fait que (∀) utilise la force de la logique classique. Par exemple, comme l'a montré J. Bell (cf. plus bas), on peut dériver de l'opérateur ε *la loi du tiers exclu* sous la forme

$$\neg \forall x\, F(x) \Rightarrow \exists x\, \neg F(x).$$

Sans l'utilisation de ¬¬F = F on ne peut déduire que *le principe de Markov*

$$(\forall x\, \neg\neg F(x) \Rightarrow F(x)) \Rightarrow (\neg \forall x\, F(x) \Rightarrow \exists x\, \neg F(x))$$

qui affirme le tiers exclu pour les prédicats *décidables* (i.e. tels que ¬¬F = F).

En ce qui concerne la sémantique des ε-termes associés aux énoncés universels, supposons que l'énoncé ∀x F(x) soit valide une fois interprété dans un certain modèle 𝔐 de 𝔏. L'ε-terme $\varepsilon_{\neg F}$ subsiste toujours symboliquement. Mais il ne peut pas dénoter de façon consistante, conformément à son sens, puisque par hypothèse il n'existe aucun individu satisfaisant ¬F. Tout référent de $\varepsilon_{\neg F}$ *nie* le contenu de l'idée dont il est l'hypostase individuée. On dit dans ce cas que $\varepsilon_{\neg F}$ est un *terme-zéro* (un "null-term"). On pourrait dire également que $\varepsilon_{\neg F}$ est en quelque sorte un symbole-index "clivé" dont l'aspect symbole (syntaxique) et l'aspect index (sémantique) sont devenus contradictoires. C'est d'ailleurs pourquoi l'interprétation de l'opérateur ε a fait problème et qu'un certain flottement était repérable chez les premiers commentateurs. Elle présuppose en effet de bien faire la différence entre, côté syntaxe, la *description-attribution* et, côté sémantique, la *référence-dénotation* (cf. Donnellan, 1966). C'est un aspect de l'opposition classique entre *Sinn* et *Bedeutung*. Un terme-zéro peut référer mais sans être aucunement attributif. La description indéfinie qu'il symbolise est alors *impropre*.

**4.    Un exemple d'ε-terme zéro : les infinitésimales leibniziennes**

Il ne faut pas croire que les termes-zéro soient des curiosités. Mon exemple favori, développé il y a longtemps dans l'*Enciclopedia Einaudi*, reste celui des *infinitésimales leibniziennes*. Soit ℝ le corps totalement ordonné des nombres réels. Sa structure d'ordre est *archimédienne* : elle satisfait l'axiome selon lequel tout nombre (aussi grand que l'on veut) est atteignable par tout nombre (aussi petit que l'on veut) à condition d'additionner ce dernier à lui-même un assez grand nombre de fois. Autrement dit, ℝ satisfait l'énoncé :

(A)     $\forall y \in \mathbb{R}^+ \, \forall x \in \mathbb{R}^+ \, \exists n \in \mathbb{N} \, (nx > y)$

(où $\mathbb{R}^+$ est l'ensemble des nombres strictement positifs et $\mathbb{N}$ l'ensemble des entiers naturels).

En passant des grands nombres à leurs inverses, cet axiome dit qu'il n'existe pas d'infinitésimale dans $\mathbb{R}$ : tout nombre non nul (positif) aussi petit que l'on veut est plus grand qu'un autre nombre non nul. Il n'existe donc pas de nombre (positif) qui soit non nul et plus petit que tous les nombres strictement positifs. Autrement dit, $\mathbb{R}$ satisfait l'énoncé :

(I)     $\forall y\ (y \neq 0) \Rightarrow \exists r\ ((r > 0) \wedge (|y| > r))$

(où $|y|$ est le symbole "valeur absolue").

Il est facile de voir que le concept leibnizien d'infinitésimale recouvre très exactement le terme-zéro associé à l'universelle (I). En effet (I) est de la forme $\forall y\ G(y)$ et est donc équivalente à l'énoncé hilbertien $G(\varepsilon_{\neg G})$. Comme $G$ est de la forme $A \Rightarrow B$, $\neg G$ est de la forme $A \wedge \neg B$ :

$\neg G(y) \equiv (y \neq 0) \wedge \forall r\ ((r > 0) \Rightarrow (|y| \leq r))$.

$\varepsilon_{\neg G}$ correspond donc à *l'idée* d'un nombre différent de 0 et dont la distance à 0 est inférieure à tout nombre réel strictement positif. C'est le *dx* leibnizien et ce *dx* est donc un terme-zéro. L'axiome d'Archimède sous la forme (I) est alors équivalent à l'énoncé :

(I')    $G(\varepsilon_{\neg G}) \equiv (\varepsilon_{\neg G} \neq 0) \Rightarrow \exists r\ ((r > 0) \wedge (|\varepsilon_{\neg G}| > r))$.

Toute infinitésimale est soit nulle soit finie : *il n'existe pas de référent numérique dans $\mathbb{R}$ pour l'idée in individuo d'infinitésimale*. Pour trouver des référents il faut passer par l'Analyse non standard (Petitot 1979).

## 5.     Le problème de la sémantique des ε-termes

L'interprétation de l'opérateur $\varepsilon$ a fait quelque peu problème à cause de ces difficultés. Un certain flottement est repérable même chez Bourbaki qui le place pourtant au fondement de sa logique :

> "Si $F$ est une assertion et $x$ une lettre, $\varepsilon_x F(x)$ est un objet; considérons l'assertion $F$ comme exprimant une propriété de l'objet $z$; alors, s'il existe un objet possédant la propriété en question, $\varepsilon_x F(x)$ représente un objet privilégié qui possède cette propriété; sinon $\varepsilon_x F(x)$ représente un objet dont on ne peut rien dire" [10].

De même chez R. Godement:

> "[L'opération de Hilbert] consiste à choisir une fois pour toutes, pour chaque relation $F$ et chaque lettre $x$, un objet vérifiant la relation $F(x)$ (s'il en "existe"; dans le cas contraire $\varepsilon_x F(x)$ est un objet dont on ne peut rien dire). Il va de soi que ce "choix" est purement fictif : l'intérêt de l'opération de Hilbert est de donner un procédé parfaitement artificiel mais purement mécanique pour construire effectivement un objet dont on sait *seulement* qu'il satisfait à des conditions imposées d'avance (dans le cas où de tels objets existeraient)." [11]

La formule ($\forall$) est à l'origine des travaux de Hilbert sur *les fonctions de choix transfinies*. Initialement, Hilbert notait $\tau_x F(x)$ l'individu $\varepsilon_{\neg F}$ et le pensait comme un *contre-exemple générique* à $F$. Les énoncés universels ($\forall$) étaient alors analogues aux universalisations fréquentes en langue naturelle : "si même lui... alors tous..." : si même $\tau_F$ qui représente idéalement un individu satisfaisant $\neg F$ satisfait quand même $F$ alors vraiment tous les individus satisfont $F$. Ensuite Hilbert changea l'interprétation de l'opérateur $\tau$ puis changea de notation et passa à l'opérateur $\varepsilon$.

L'opérateur de Hilbert rend ainsi dans une certaine mesure la contradiction opératoire, non pas la contradiction logique brute mais cette forme plus subtile de contradiction qu'est l'incompatibilité entre le signifié du symbole-index $\varepsilon_F$ et les propriétés "réelles" de ses référents. Il conduit à prendre prendre en compte des distinctions inhabituelles en théorie des modèles.

(i)     Dans un individu-type $\varepsilon_F$ il faut distinguer — comme pour les énoncés — la structure syntaxique (la construction de $\varepsilon_F$ à partir de $F$), le sens (ou le signifié, le contenu idéal, le *Sinn*) et les référents (la dénotation, la *Bedeutung*).

(ii)    Syntaxiquement (symboliquement), $\varepsilon_F$ subsiste. Lorsqu'il admet des référents conformes à son sens, alors il existe. Sinon sa subsistance n'est qu'une existence purement symbolique.

(iii)   Il existe une double nature des $\varepsilon$-terme. $\varepsilon_F$ symbolise un "ceci" indéterminé "qui satisfait $F$", un "celui qui...". Mais il réfère aussi à un élément, au demeurant quelconque, de l'extension $X_F$ si $X_F \neq \emptyset$.

C'est dans ce jeu entre la négation, l'article indéfini et le démonstratif que se joue l'un des principaux intérêts de l'opérateur de Hilbert (Klaus von Heusinger y a récemment aussi beaucoup insisté). Et cela d'autant plus que participe également au jeu l'article *défini* dans son usage *non anaphorique*. Il existe en effet une troisième interprétation de l'opérateur $\varepsilon$ — l'interprétation que j'ai proposé d'appeler *générique* — qui consiste à interpréter $\varepsilon_F$ comme dénotant un élément générique de l'extension $X_F$ (si $X_F \neq \emptyset$). La relation primitive d'appartenance $a \in X_F$ se trouve alors substituée par une relation (également primitive) de *spécialisation* $\varepsilon_F \to a$ : l'élément $a \in X_F$ est un individu qui spécifie et particularise — spécialise — l'élément générique $\varepsilon_F$. Dans cette interprétation, si un prédicat $G(y)$ est valide pour l'élément générique $\varepsilon_F$, il est *ipso-facto* valide pour toute ses spécialisations, autrement dit, on a les équivalences :

$$G(\varepsilon_F) \Leftrightarrow \varepsilon_F \in X_G \Leftrightarrow \forall x\,[(\varepsilon_F \to x) \Rightarrow x \in X_G] \Leftrightarrow X_F \subseteq X_G \Leftrightarrow \forall x\,[F(x) \Rightarrow G(x)].$$

---

[10] Bourbaki [1958].

[11] Godement [1962].

On voit ainsi en définitive que, comme beaucoup d'autres entités intensionnelles, les $\varepsilon$-termes possèdent une interprétation *de dicto* (l'interprétation générique) et une interprétation *de re* (l'interprétation spécifique comme fonction de choix).

Ces divers aspects de l'opérateur $\varepsilon$ deviennent particulièrement intéressants et significatifs si l'on songe qu'à partir de lui Hilbert et ses collaborateurs ont élaboré un calcul logique — dit $\varepsilon$-calculus et noté $CP_\varepsilon$ — qui s'est révélé être *syntaxiquement équivalent au calcul classique des prédicats* CP. $CP_\varepsilon$ est constitué des composantes suivantes :

(i)     le calcul standard des propositions;
(ii)    des symboles de constantes individuelles et des symboles de variables en nombre suffisant;
(iii)   des symboles de prédicats *n*-aires, pour tout *n*;
(iv)    les $\varepsilon$-termes correspondants;
(v)     les règles standard de déduction dans CP;
(vi)    la règle d'introduction de l'opérateur $\varepsilon$ ($\varepsilon$-formula) :

($\varepsilon$)     $F(a) \Rightarrow F(\varepsilon_F)$.

On démontre alors (cf. plus bas) que si $\varphi$ est un énoncé de CP et $\psi$ une conséquence de $\varphi$ dans $CP_\varepsilon$ d'où le symbole $\varepsilon$ est éliminable, alors $\psi$ est en fait une conséquence de $\varphi$ dans CP. Autrement dit, $CP_\varepsilon$ est une extension *inessentielle* (on dit aussi *conservative*) de CP. Ce résultat n'est plus vrai du tout en logique intuitionniste.

Il n'y a donc aucune raison de considérer CP comme plus évident ou plus naturel que $CP_\varepsilon$. En fait, $CP_\varepsilon$ est même philosophiquement beaucoup plus intéressant car, tout en étant conservatif sur CP, il est plus expressif que ce dernier dans la mesure où les $\varepsilon$-termes peuvent être combinés de façon plus compliquée que les quantificateurs.

### III. LES ORIGINES METAMATHEMATIQUES ET LA PORTEE EPISTEMOLOGIQUE DE L'OPERATEUR $\varepsilon$

#### 1. La stratégie finitiste hilbertienne

L'opérateur $\varepsilon$ a été introduit par Hilbert comme un outil dans son vaste programme formaliste qui consistait à démontrer par des moyens *finitistes* la *consistance* des principales théories mathématiques et, avant tout, de l'arithmétique formelle. On sait que ce programme de recherche a dominé la logique mathématique jusqu'à ce que le théorème d'incomplétude de Gödel en ait démontré les limites intrinsèques. Ce résultat négatif a en partie fait déchoir pour les mathématiciens les tentatives pré-gödeliennes au rang de curiosités, mais cela ne les empêche pas de garder tout leur intérêt cognitif, sémiologique et philosophique.

On sait que la thèse de Hilbert était que, pour "fonder" les mathématiques, il fallait démontrer leur consistance en ne faisant usage que d'une métalogique *finie*, celle-ci étant la seule à pouvoir être considérée comme légitime et évidente *a priori*. Mais cela soulevait immédiatement le problème de la présence de *quantificateurs* dans les axiomes de théories possédant des modèles *infinis*. Car alors la quantification n'est plus un processus finitiste: $\exists x\, F(x)$ (resp. $\forall x\, F(x)$) n'est plus une simple disjonction finie $F(a_1)\vee...\vee F(a_n)$ (resp. une simple conjonction finie $F(a_1)\wedge...\wedge F(a_n)$. Il fallait donc pouvoir *éliminer* les quantificateurs des axiomes sans pour autant affaiblir la force du calcul des prédicats utilisé. C'est pourquoi Hilbert a cherché à *définir* les quantificateurs dans le cadre d'un formalisme satisfaisant aux conditions très restrictives imposées par sa stratégie finitiste. Et il a cru trouver le cadre approprié avec l'$\varepsilon$-calcul $CP_\varepsilon$ (et ses dérivés). En effet, la stratégie finitiste impose que l'on se restreigne à des manipulations logiques élémentaires sur des énoncés *élémentaires* (sans quantificateurs) de type $F(a)$, $G(a,b)$, etc. C'est donc en *élargissant* l'univers des objets par l'introduction d'objets idéaux qu'il faut arriver à éliminer les quantificateurs. Ce qui est possible dans $CP_\varepsilon$ avec l'introduction de la "transfinite logische Auswahlfunktion" qu'est l'opérateur $\varepsilon$.

Après quelques essais remontant à 1923, Hilbert a introduit sa fonction de choix transfinie dans son célèbre mémoire de 1925 *Sur l'Infini* [12] et, avec Bernays, il en a développé l'usage dans son ouvrage monumental *Grundlagen der Mathematik*. Il a accédé grâce à lui à un certain nombre de démonstrations constructives et finitistes de consistance (cf. plus bas). Leisenring décrit ainsi sa stratégie formaliste générale.[13]

---

[12] Hilbert [1925]. Cf. aussi Hilbert [1927] ainsi que Bernays [1927] et Ackermann [1928].

> "Supposons que 𝔗 soit une théorie fondée sur le calcul des prédicats. S'il existe un modèle 𝔐 qui satisfait l'ensemble 𝔄 des axiomes de 𝔗, alors 𝔗 est consistante. […] L'objection à ce type de preuve de consistance est qu'elle requiert une très forte métathéorie. Par exemple, si la cardinalité du modèle est infinie, il faut, pour démontrer que $\phi$ [une formule fausse quelconque, par exemple $0 = 1$] n'est pas un théorème de 𝔗 [et que 𝔗 est donc consistante], utiliser des arguments non-constructifs montrant que les axiomes de CP sont vrais dans le modèle. Une des principales contributions des formalistes a été de montrer que pour certaines théories particulières ce type de preuve de consistance pouvait être menée à bien de façon complètement finitiste."[14]

Considérons par exemple les axiomes de la théorie 𝔗 des ensembles infinis sur lesquels est définie une relation < d'ordre total avec un plus petit élément.

$A_1$ $\forall x \neg (x < x)$ non réflexivité de l'ordre
$A_2$ $\forall x \forall y \forall z ((x < y) \wedge (y < z) \Rightarrow (x < z))$ transitivité de l'ordre
$A_3$ $\forall x \forall y ((x < y) \vee (y < x) \vee (x = y))$ ordre total
$A_4$ $\forall x \exists y (x < y)$ ordre infini
$A_5$ $\exists x \forall y ((x = y) \vee (x < y))$ plus petit élément

L'ensemble ℕ des nombres entiers muni de sa relation d'ordre naturel est un modèle de 𝔗 et donc 𝔗 est consistante. Mais cette preuve sémantique est non constructive. Pour accéder à une preuve syntaxique finitiste, il faut avant tout éliminer le quantificateur ∃ des axiomes $A_4$ et $A_5$. On peut le faire par la méthode de "résolution" symbolique suivante. Soit un axiome de la forme $\exists x \forall y\, B(x,y)$. On peut lui substituer $\forall y\, B(s,y)$ (où $s = \varepsilon_x \forall y\, B(x,y)$ est un $\varepsilon$-terme) qui est un axiome ne contenant plus qu'un quantificateur universel. De même, soit un axiome de la forme $\forall x \exists y\, B(x,y)$. Suivant Skolem, introduisons un symbole de fonction $g$ (dite de Skolem) et posons $g(x) = \varepsilon_y B(x,y)$. On a alors :

$$\forall x \exists y\, B(x,y) \equiv \forall x\, B(x,\varepsilon_y B(x,y)) \equiv \forall x\, B(x,g(x)).$$

Dans notre cas on obtient donc le système d'axiomes :

---

[13] Leisenring [1969]. Pour une introduction à la théorie des modèles, cf. Petitot [1979].

[14] Leisenring [1969], pp.85-86.

$A_1, A_2, A_3, A'_4 : \forall x\, (x < g(x))$, $A'_5 : \forall y\, ((s = y) \vee (s < y))$,

axiomes d'où les quantificateurs $\exists$ ont été éliminés. Ces axiomes sont vrais dans $\mathbb{N}$ en prenant pour $g$ la fonction successeur $g(x) = x+1$ et pour $s$ le plus petit élément $s = 0$.

On utilise alors le résultat suivant. On sait que toute formule $\varphi$ de CP peut être mise sous une forme dite "prénexe" qui la rend équivalente à une formule $F'$ de la forme $\Pi\varphi^0$ où $\Pi$ est une suite de quantificateurs et $\varphi^0$ une formule *sans* quantificateurs dite "matrice" de $\varphi$. *Le premier théorème fondamental d'élimination de $\varepsilon$* (cf. plus bas) dit que si $\Sigma$ est un ensemble de formules prénexes et si $\varphi$ est une formule prénexe, alors si $\Sigma \vdash_{CP} \varphi$ (i.e. si $\varphi$ est une conséquence de $\Sigma$ dans CP) on a $\Sigma^* \vdash_{CE} \varphi_1 \vee \ldots \vee \varphi_n$ où les formules de $\Sigma^*$ (resp. les formules $\varphi_i$) sont obtenues à partir des matrices des formules de $\Sigma$ (resp. de la matrice de $\varphi$) en y substituant les variables par des termes (y compris par des $\varepsilon$-termes) et où CE est le calcul *élémentaire* des prédicats ne faisant plus usage des quantificateurs. Soit alors $\Sigma$ le système d'axiomes $A_1, \ldots, A'_5$. Les matrices en sont :

$$\begin{cases} A_1^0 & \neg(x < x) \\ A_2^0 & (x < y) \wedge (y < z) \Rightarrow (x < z) \\ A_3^0 & (x < y) \vee (y < x) \vee (x = y) \\ A_4^0 & x < g(x) \\ A_5^0 & (s = y) \vee (s < y)) \end{cases}$$

Pour démontrer *syntaxiquement* que $\mathfrak{T}$ est consistante, il suffit donc d'assigner une interprétation aux symboles et une valeur de vérité aux formules *élémentaires* du langage formel considéré qui soient telles que toute substitution dans une de ces matrices donne un énoncé vrai. Or cela est facile à effectuer de façon finitiste. Comme l'explique Leisenring :

> "De façon générale, la méthode formaliste des preuves de consistance peut être décrite de la façon suivante. Supposons que $\mathfrak{T}$ soit une théorie basée sur CP. En remplaçant chaque axiome de $\mathfrak{T}$ par une forme prénexe équivalente, en prenant les résolutions de Skolem de ces formules prénexes et en adjoignant les nouvelles fonctions de Skolem au vocabulaire [du langage formel], on obtient une théorie $\mathfrak{T}'$ qui est une extension inessentielle de $\mathfrak{T}$. On essaye alors de trouver une assignation effective de valeurs de vérités aux formules atomiques de $\mathfrak{T}'$ d'une façon telle que chaque axiome $E_1$ et chaque axiome $E_2$ [les axiomes de l'égalité dans CP] et que chaque substitution des matrices des axiomes de $\mathfrak{T}'$ prennent la valeur 1 [i.e. soient vrais]. Si cela peut être fait alors à la fois

𝔗´ et 𝔗 sont consistantes."[15]

C'est de cette façon (et en particulier en définissant explicitement à partir de l'opérateur $\varepsilon$ les fonctions de Skolem) que Hilbert a justifié l'usage de l'infini. Comme le note encore Leisenring :

> "Peut-être que la signification majeure du premier $\varepsilon$-théorème de Hilbert est la suivante. Bien que la logique classique, telle qu'elle est formalisée par le calcul des prédicats, contient des aspects non-finitistes, toute preuve d'un énoncé finitiste peut être convertie en une preuve finitiste."[16]

Toutefois, comme nous le verrons plus bas, il s'est heurté à une difficulté irréductible qui n'a été vraiment éclaircie que par Gödel.

**2. L'interprétation philosophique d'Albert Lautman**

Cette stratégie métamathématique a été admirablement décrite sur le plan épistémologique par Albert Lautman.[17] Dans plusieurs textes, ce dernier est revenu sur le fait que le formalisme hilbertien échappe aussi bien au logicisme russellien qu'aux contraintes trop drastiques de l'intuitionnisme brouwerien.

Le problème est, nous l'avons vu, que pour les modèles infinis de théories il y a contradiction entre l'exigence de manipulation finitiste des symboles (des variables) intervenant dans les formules et le nombre infini de vérifications élémentaires qu'exige la validation de formules contenant des quantificateurs. D'où l'idée d'introduire

> "des champs métamathématiques, intermédiaires entre les signes des formules et les champs mathématiques [les modèles] de leurs valeurs effectives."[18]

A la suite de Hilbert, Herbrand cherche à ramener tout énoncé à un énoncé sans quantificateurs dont la valeur de vérité soit vérifiable en un nombre fini d'étapes. Pour cela

> "il lui a fallu trouver une manière de définir, pour une variable susceptible de prendre une infinité de valeurs mathématiques, un nombre fini de valeurs métamathématiques, qui symbolisent ainsi l'existence de

---

[15] *Ibid.*, p.87.

[16] *Ibid.*, p.88.

[17] Pour une introduction à la philosophie mathématique de Lautman, cf. Petitot [1987] et [2001].

[18] Lautman [1937], p. 108.

cette infinité si difficile à manier."[19]

C'est là que les procédures de type opérateur $\varepsilon$ deviennent essentielles. Si elles permettent d'élaborer des preuves constructives (syntaxiques et effectives) qui échappent cependant aux contraintes trop strictes des constructivismes intuitionnistes, c'est parce qu'elles imposent les contraintes de finitude aux entités *métamathématiques* et non pas directement aux entités des modèles considérés. Ces entités métamathématiques que sont les $\varepsilon$-termes, les fonctions de Skolem, etc. sont, selon Lautman, des "*mixtes*" intermédiaires entre les symboles et les référents dénotés par ces symboles, c'est-à-dire *entre syntaxe et sémantique*.

> "Les éléments de ces champs [métamathématiques] sont donc en étroite correspondance avec les signes des variables des formules; ils constituent plutôt un système de nouveaux signes que l'on substitue aux premiers qu'un ensemble de véritables valeurs pour les variables désignées par ces signes. D'un autre côté, ils n'en possèdent pas moins une nature de champs indépendants de la formule qu'ils réalisent, et présentent bien ainsi un premier aspect de mixtes en tant qu'ils sont intermédiaires entre les signes formels et leurs valeurs mathématiques effectives. (...) Intermédiaires entre les signes et leurs valeurs, ces champs sont, d'une part homogènes à la discontinuité finie des signes puisqu'à un signe de variable ($\exists x$) ne correspond qu'une valeur $a$ [$a$ est ici l'$\varepsilon$-terme $\varepsilon_F$ associé à la formule $F(x)$ que $\exists x$ quantifie] et, d'autre part ils symbolisent une infinité de valeurs mathématiques puisque la lettre $a$ représente n'importe quelle valeur mathématique éventuelle de la variable $x$ lorsqu'elle intervient sous la forme particulière ($\exists x$). Une médiation s'opère donc par ces champs du fini à l'infini, qui permet (...) de dominer l'infini."[20]

Dans son dernier essai, *Considérations sur la logique mathématique*, Lautman est revenu sur ce point décisif.

> "Il est impossible d'opérer tous les calculs impliqués par une théorie [à modèle infini], car ils sont évidemment en nombre infini et les intuitionnistes ont raison de dire qu'on ne serait jamais certain en procédant ainsi de ne pas trouver de contradiction. Mais il est possible de remplacer, en ce qui concerne l'étude de la valeur logique, la considération d'une infinité de valeurs particulières par une lettre ou une fonction "*de choix*" tels que les résultats obtenus dans le champ fini de ces valeurs de choix aient des conséquences valables transfiniment pour tous les êtres mathématiques particuliers dont les valeurs sont symbolisées par cette valeur de choix. (...) Entre les exigences de construction des intuitionnistes et la pure introduction de notions par axiomes, les hilbertiens ont su interpréter un schématisme intermédiaire, celui d'individu et de champs considérés non tant pour eux-mêmes que pour les conséquences infinies que permettent les calculs finis opérés grâce à eux."[21]

---

[19] *Ibid*.

[20] *Ibid*., p.109.

### 3.   ε-calculus et objets idéaux

Avant de préciser tous ces points, faisons une remarque générale. L'$\varepsilon$-calcul hilbertien est un calcul des prédicats où l'univers des objets est élargi par l'adjonction d'objets idéaux génériques. Dans $CP_\varepsilon$ ces objets sont canoniquement associés aux formules, mais l'on peut, de façon plus générale, envisager d'autres types d'adjonction. Sur ce point, on pourra consulter par exemple, les travaux de Kit Fine.[22] Kit Fine insiste en particulier sur les faits suivants.

(i)   Une interprétation en terme d'objets génériques permet immédiatement de trouver les restrictions qui doivent être apportées, dans les divers systèmes logiques, aux règles d'instanciation (I) et de généralisation (G) respectivement universelle (U) et existentielle (E) :

$$\text{IU} \quad \frac{\forall x F(x)}{F(a)} \qquad \text{GU} \quad \frac{F(a)}{\forall x F(x)}$$

$$\text{IE} \quad \frac{\exists x F(x)}{F(a)} \qquad \text{GE} \quad \frac{F(a)}{\exists x F(x)}$$

Les règles IU et GE étant évidentes, ce sont évidemment les règles GU et IE qui font problème. Par exemple, comme nous l'avons déjà noté plus haut, à cause de l'existence de modèles non standard, on peut parfaitement pouvoir démontrer pour chaque $a$ l'énoncé $F(a)$ sans pouvoir pour autant démontrer l'universelle $\forall x\, F(x)$ (cf. Petitot [1995] pour l'exemple de la version finitiste du théorème de Kruskal).

Par exemple, dans le système Q proposé par Quine on a les restrictions suivantes :
$R_1$ :   le symbole d'objet (générique) $a$ n'intervient pas dans $F(x)$ dans les règles GU et IE. Dans une application de GU ou de IE, soient $b$, $c$, etc. les constantes de $F$ et soit $a$ ($a \neq b$, $c$, etc.) une instanciation de $x$. On dira que $a$, terme instancié, dépend immédiatement des termes donnés $b$, $c$, etc.
$R_2$ :   un symbole $a$ d'objet générique ne peut pas être instancié deux fois.

---

[21] *Ibid.*, pp.313-314.

[22] Cf. Fine [1985]. Un excellent exposé de ces questions se trouve aussi dans Nef [1998].

R$_3$ : on peut ordonner les termes instanciés dans une démonstration dans un ordre $a_1, ..., a_n$ tel que pour tout $i$ aucun des $a_{i+1}, ..., a_n$ ne dépende immédiatement de $a_i$.

Ces restrictions sont naturelles si l'on interprète dans Q, le $a$ de IE comme un représentant typique et générique de $X_F$ (i.e. comme $\varepsilon_F$, IE devenant l'équivalence ($\exists$) de CP$_\varepsilon$) et le $a$ de GU comme un "contre-exemple" générique de $F$ (i.e. comme $\tau_F = \varepsilon_{\neg F}$, GU devenant l'équivalence ($\forall$) de CP$_\varepsilon$).

(ii) Ces interprétations en termes de sémantique de la généricité permettent de comprendre les différences importantes qui existent dans les divers systèmes logiques entre les restrictions apportées aux règles ci-dessus. Par exemple une autre interprétation de GU consiste à poser que $a$ est l'objet générique de tout l'univers d'objet considéré (par exemple $a = \varepsilon_x(x = x)$). Cela est très différent de l'interprétation quinienne de $a$ comme $\tau_F = \varepsilon_{\neg F}$.

(iii) Les objets génériques permettent de donner un sens plus clair aux fonctions de Skolem. Par exemple, une instanciation de IE $\dfrac{\exists x F(x, b)}{F(a, b)}$ sera représentée par $\dfrac{\exists x F(x, b)}{F(f(b), b)}$ où $f$ est la fonction de Skolem associée. On sait que les avis divergent sur l'interprétation de $f$. Est-ce une fonction définie ($f$ étant alors un symbole de constante) ou indéfinie? Qu'elle soit indéfinie est contraire à l'intuition. Mais en général, elle ne peut pas être explicitement définie. Elle le devient toutefois si sa valeur $f(b)$ est considérée comme un terme générique.

(iv) Enfin, et c'est un point décisif, la pensée et le raisonnement *naturels* reposent sur l'usage d'objets génériques. La logique naturelle est une logique *catégorique de la généricité* et non pas une logique ensembliste extensionnelle. Comme l'affirme Kit Fine

> "le principal avantage de la sémantique générique est qu'elle fournit une représentaiton fidèle du raisonnement ordinaire."[23]

## IV. L'$\varepsilon$-CALCULUS

### 1. Principes du calcul

Donnons maintenant quelques précisions un peu plus techniques sur les sources et les développements logico-mathématiques de l'idée hilbertienne.

L'affaire remonte à Russell et à son opérateur $\iota$ formalisant la notion de *description définie*. Considérons une description définie, autrement dit une formule $F(x)$ qui caractérise un individu parce qu'elle satisfait aux deux conditions d'existence et d'unicité :

---

[23] Fine [1985].

(E)   $\exists x\, F(x)$

(U)   $\forall x \forall y\, (F(x) \wedge F(y) \Rightarrow x = y)$.

Russell introduit alors la notation $\iota_x F(x)$ ($= \iota_F$) pour symboliser l'individu caractérisé par $F$ et élargit CP en lui adjoignant la règle :

(R$\iota$) :   sous les conditions (E) et (U) on peut introduire dans une démonstration l'individu $\iota_F$ et inférer $F(\iota_F)$.

     L'opérateur $\iota$ ne pose pas de problème d'interprétation. Il formalise *l'article défini*. Et comme $\#X_F = 1$, l'article défini (*le* $a \in X_F$), l'article indéfini (*un* $a \in X_F$ quelconque) et le démonstratif (*celui qui* satisfait $F$) coïncident. De par sa structure c'est bien un symbole-index, mais, puisque $\#X_F = 1$, sa dénotation ne comporte aucune ambiguïté.

     Lorsque son emploi est licite, l'opérateur $\iota$ permet immédiatement d'effectuer les procédures de résolution symbolique par introduction de fonctions de Skolem (élimination des quantificateurs, cf. plus haut). Soit par exemple un axiome $A$ de la forme $A \equiv \forall x \exists y\, B(x,y)$ et supposons que dans le système considéré l'on puisse dériver la propriété d'unicité : $\forall x \forall y \forall z\, (B(x,y) \wedge B(x,z) \Rightarrow y = z)$. Par application de (R$\iota$) on peut introduire $\iota_y B(x,y)$ quel que soit $x$ et, par introduction d'un symbole de fonction $g$, poser $g(x) = \iota_y B(x,y)$. On aura alors, comme plus haut, $A \equiv \forall x\, B(x, g(x))$, la fonction de Skolem $g$ étant ici *explicitement définie*.

     La première généralisation effectuée par Hilbert a consisté à s'affranchir de l'hypothèse d'unicité (U). Sous la seule condition d'existence (E), Hilbert introduit un symbole $\eta$ et un $\eta$-terme $\eta_x F(x)$ dont l'introduction est régie par la règle :

(R$\eta$) :   sous la condition (E) on peut introduire dans une démonstration l'individu $\eta_F$ et inférer $F(\eta_F)$.

Avec les $\eta$-termes s'introduisent dans CP, selon le dire même de Hilbert, des éléments idéaux (des idées *in individuo*) canoniquement associés à des formules. L'opérateur $\eta$ formalise en particulier *l'article indéfini*. Sa dénotation consiste en un élément, au demeurant quelconque, de $X_F$ (non vide par hypothèse). $\eta_F$ est déjà un "vrai" symbole-index avec la triple interprétation : article indéfini, article défini non anaphorique (interprétation comme objet générique, i.e. *de dicto*), démonstratif ("celui qui satisfait $F$" : interprétation *de re*). Hilbert a privilégié cette troisième interprétation en considérant qu'un $\eta$-terme $\eta_F$ était un *opérateur de choix*.

Hilbert a alors subtilement remarqué que, indépendamment de toute condition d'existence et d'unicité, l'énoncé : $\exists x\,(\exists y\,F(y) \Rightarrow F(x))$ est *toujours* dérivable. On peut donc *toujours* introduire l'$\eta$-terme $\varepsilon_x F(x)$ (correspondant à la règle IE) :

(*) $\quad \varepsilon_x F(x) \equiv \eta_x(\exists y\,F(y) \Rightarrow F(x))$.

Si la condition d'existence (E) est satisfaite, alors $\varepsilon_F = \eta_F$. Si de plus, la condition d'unicité (U) est satisfaite, alors $\varepsilon_F = \eta_F = \iota_F$. L'opérateur $\varepsilon$ est donc bien une généralisation des descriptions définies russelliennes. Mais les $\varepsilon$-termes sont *des individus idéaux qui involuent dans leur symbole même la question de leur existence*. C'est ce "tour" éminemment métaphysique qui, selon nous, donne à l'opérateur de Hilbert sa portée philosophique.

$\varepsilon_F$ est par définition un $\eta$-terme $\eta_x G(x)$ avec $G(x) \equiv \exists y\,F(y) \Rightarrow F(x)$. D'après la règle ($R_\eta$) on peut donc *toujours* dériver l'énoncé $G(\eta_G)$, soit $\exists y\,F(y) \Rightarrow F(\varepsilon_F)$. Comme il est clair que, si $F(\varepsilon_F)$ est vrai alors, par généralisation existentielle GE, on a $\exists y\,F(y)$, on a bien en définitive *l'équivalence logique* :

($\exists$) $\quad \exists y\,F(y) \Leftrightarrow F(\varepsilon_F)$.

C'est ce constat qui a conduit Hilbert à inverser la démarche, à introduire d'emblée l'opérateur $\varepsilon$ comme une *primitive* du calcul logique et, comme nous l'avons vu, à définir la quantification existentielle par ($\exists$). Il a ainsi obtenu sa variante $CP_\varepsilon$ du calcul des prédicats. $CP_\varepsilon$ contient un axiome spécifique remplaçant la règle GE d'introduction du quantificateur $\exists$ dans CP : de $F(a)$ (où $a$ n'apparaît pas dans $F$) on peut inférer $\exists x\,F(x)$. Cet axiome est l'*$\varepsilon$-formula* évoquée plus haut :

($\varepsilon$) $\quad F(a) \Rightarrow F(\varepsilon_F)$.

On lui adjoint éventuellement *l'axiome d'égalité*, dû à Ackermann :

($\varepsilon_2$) $\quad \forall x\,(F(x) \Leftrightarrow G(x)) \Rightarrow \varepsilon_F = \varepsilon_G$.

Il est essentiel de noter que ce dernier axiome formule un principe *d'extensionalité* qui est très loin d'être évident. Son évidence dépend de l'interprétation de $\varepsilon$. Il signifie que si $F$ et $G$ sont extensionnellement équivalentes, i.e. si au niveau de leurs extensions $X_F = X_G$, alors $\varepsilon_F = \varepsilon_G$. Cela n'est déjà pas évident dans l'interprétation générique *de dicto* car comme les sens de $F$ et de $G$ sont en général différents on ne voit pas pourquoi les éléments génériques satisfaisant typiquement $F$ et $G$ devraient être identiques. On devrait même poser plutôt $\varepsilon_F \neq \varepsilon_G$ de façon à rendre compte de la différence de sens (*Sinn*) de $F$ et de $G$ même si leur dénotations (*Bedeutung*) sont équivalentes. Cela est encore moins évident dans l'interprétation spécifique *de re*. On ne voit pas pourquoi en effet, les objets sélectionnés par $\varepsilon_F$ et $\varepsilon_G$ devraient nécessairement être identiques. L'intérêt principal des ε-termes est d'être des entités intensionnelles et même plus précisément, nous l'avons vu plus haut, des objets intensionnellement décomplétés. L'axiome d'Ackermann en refait subrepticement des entités extensionnelles. Nous reviendrons en conclusion sur ce point crucial.

Nous retrouvons ainsi les problèmes d'interprétation sémantique de l'opérateur $\varepsilon$. Jusqu'à une date relativement récente, on n'a pas tenu compte de sa nature intensionnelle et on l'a essentiellement interprété, nous l'avons vu, en termes de fonctions de choix. Soit $\mathcal{U}$ l'univers du discours (le modèle) considéré. On introduit *une fonction de choix universelle* $\Phi$ qui associe à chaque ensemble $X \neq \emptyset$ de $\mathcal{U}$ un élément $\Phi(X)$ qui appartient à $X$. Si $F(x)$ est un prédicat d'extension $X_F \neq \emptyset$, on pose a priori que $\varepsilon_F$ dénote $\Phi(X_F)$. Mais alors comment définir $\Phi(\emptyset)$. Selon Asser et Hermès, il faut poser $\Phi(\emptyset) = \Phi(\mathcal{U})$.[24] Revenons en effet à la définition de $\varepsilon$ à partir de $\eta$ :

(*)   $\varepsilon_x F(x) \equiv \eta_x(\exists y\, F(y) \Rightarrow F(x)) \equiv \eta_x G(x).$

Si $F$ est d'extension vide ($X_F = \emptyset$), alors pour tout $x$ l'implication $\exists y\, F(y) \Rightarrow F(x)$ est vérifiée puisque $\exists y\, F(y)$ est faux. Donc $\varepsilon_x F(x) = \eta_x G(x)$ avec $X_G = \mathcal{U}$. Il est par conséquent "naturel" de poser $\Phi(\emptyset) = \Phi(\mathcal{U})$. En particulier, on est ainsi conduit à identifier entre eux tous les termes-zéro, par exemple à "l'objet impossible" $\varepsilon_x(x \neq x)$ et à les identifier tous à "l'objet en général" $\varepsilon_x(x = x)$.

On voit à quel point dans le calcul hilbertien la consistance logique devient compatible avec l'inconsistance ontologique. Il suffit de songer au fait que, chez Frege, le concept d'un objet différent de lui-même ("l'objet impossible") est identifié à 0, et que le concept d'un objet identique à lui-même ("l'objet en général") est identifié à 1. Dans $CP_\varepsilon$ on rencontre ainsi une sorte de $0 = 1$ conceptuel qui ne comporte aucune contradiction logique (ce qui montre d'ailleurs de façon spectaculaire les limites de la conception frégéenne du concept).

Quoi qu'il en soit, muni des deux axiomes ($\varepsilon$) et ($\varepsilon_2$) et de cette interprétation, l'$\varepsilon$-

calcul $CP_\varepsilon$ est, comme CP, *adéquat et complet*, autrement dit $\Sigma \vdash_{CP_\varepsilon} \varphi$ si et seulement si $\Sigma \vDash \varphi$ (i.e. ssi $\varphi$ est valide dans tout modèle de $\Sigma$).

Venons-en donc maintenant aux théorèmes fondamentaux qui légitiment l'emploi de l'opérateur de Hilbert. Nous avons déjà énoncé plus haut le *premier $\varepsilon$-théorème*. Il dit essentiellement que s'il existe une démonstration $\Sigma \vdash_{CP} \varphi$ dans CP ou $CP_\varepsilon$ d'une formule $\varphi$ de CE à partir d'axiomes $\Sigma$ de CE, alors il en existe déjà une dérivation *dans CE*. Il est évidemment essentiel à la stratégie finitiste hilbertienne car il permet de ramener des démonstrations de consistance à des procédures finitistes de vérification. Supposons en effet, comme plus haut, que les substitutions des matrices des axiomes d'une théorie $\mathfrak{T}$ soient vérifiées (valides) dans une certaine interprétation du langage $\mathfrak{L}$ de $\mathfrak{T}$ et une certaine assignation de valeurs de vérités aux formules élémentaires de $\mathfrak{L}$. Alors tous les énoncés de CE qui en sont dérivables dans CP ou $CP_\varepsilon$ sont valides puisqu'ils sont dérivables en fait dans CE et que, clairement, les inférences dans CE préservent la validité.

Quant au *deuxième $\varepsilon$-théorème*, il affirme essentiellement que, bien qu'apparemment plus puissant que CP, $CP_\varepsilon$ est en fait — en logique *classique* — une *extension inessentielle* (conservative) de CP : soit $\Sigma$ un ensemble de formules de CP et $\varphi$ une formule sans symbole $\varepsilon$. Si $\Sigma \vdash_{CP_\varepsilon} \varphi$ alors $\Sigma \vdash_{CP} \varphi$.

Les deux $\varepsilon$-théorèmes permettent de simplifier considérablement (sans toutefois réaliser le "rêve" de Hilbert) les démonstrations syntaxiques de consistance. En effet, ils impliquent que, lorsque l'on passe d'un système d'axiomes $\Sigma$ au système $\Sigma'$ obtenu comme plus haut par résolution symbolique et élimination des quantificateurs, on obtient en fait une extension inessentielle. Soit en effet $\varphi$ une formule sans symbole $\varepsilon$ dérivable de $\Sigma'$ dans $CP_\varepsilon$ : $\Sigma' \vdash_{CP_\varepsilon} \varphi$. Comme $\Sigma \vdash_{CP_\varepsilon} \Sigma'$, on a $\Sigma \vdash_{CP_\varepsilon} \varphi$. Mais comme $\Sigma$ appartient à CP, on a d'après le deuxième $\varepsilon$-théorème $\Sigma \vdash_{CP} \varphi$. Or, d'après le premier $\varepsilon$-théorème, la démonstration de la consistance de $\Sigma'$ peut partiellement se ramener à une démonstration dans CE.

Il existe une preuve sémantique non constructive triviale du deuxième $\varepsilon$-théorème. Soit en effet $\Sigma \vdash_{CP_\varepsilon} \varphi$. Il est facile de montrer que l'$\varepsilon$-calculus $CP_\varepsilon$ est "sain" ou "adéquat", i.e. que la dérivation y conserve la validité. Donc $\Sigma \vDash \varphi$. Mais le calcul des prédicats CP étant *complet*, on a $\Sigma \vdash_{CP} \varphi$.

Mais si l'on veut obtenir une preuve syntaxique, on se heurte à la difficulté suivante. $CP_\varepsilon$ s'obtient à partir de CP en lui adjoignant les schémas d'axiomes ($\varepsilon$) et ($\varepsilon_2$). Or, si l'on substitue dans ces axiomes les $\varepsilon$-termes par des symboles de constantes, les formules obtenues *ne sont plus* des axiomes en général. Considérons par exemple une occurrence du schéma d'égalité extensionnel d'Ackermann ($\varepsilon_2$) et remplaçons-y $\varepsilon_G$ par $a$. La formule

$$\forall x \, (F(x) \Leftrightarrow G(x)) \Rightarrow \varepsilon_F = a$$

est trivialement fausse en général. Il s'agit-là, en quelque sorte, d'une "revanche de la référence". Un phénomène analogue peut se produire pour une tautologie comprenant des $\varepsilon$-termes. Considérons par exemple la tautologie $\neg \exists x \, F(x) \Rightarrow \neg F(t)$ (où $t$ est un terme quelconque) et son occurrence pour $F(x) \equiv (x = \varepsilon_y(y = x))$. On obtient :

$$\neg \exists x \, (x = \varepsilon_y(y = x)) \Rightarrow \neg(t = a) \Leftrightarrow t \neq a$$

qui est fausse en général. C'est ce phénomène, dit *d'impropreté*, qui rend le deuxième $\varepsilon$-théorème non trivial.

---

[24] Cf. Asser [1957].

Cependant soit $CP_\varepsilon^*$ l'$\varepsilon$-calcul *restreint aux formules propres*. Il est alors facile d'y démontrer syntaxiquement le deuxième $\varepsilon$-théorème.[25] Soit en effet $\varphi$ une formule sans symbole $\varepsilon$ dérivable dans $CP_\varepsilon^*$. Soit $\Delta$ une démonstration de $\varphi$ et $n$ le nombre d'$\varepsilon$-termes intervenant dans $\Delta$. On raisonne par récurrence. Soit $\varepsilon_y B(y)$ dans $\Delta$, $B$ étant sans $\varepsilon$-termes. Soit $a$ un symbole de constante n'intervenant ni dans $B$ ni dans $\varphi$. Soit $\Delta'$ la séquence obtenue par la substitution $\varepsilon_y B(y) \to a$. Si $\exists y \, B(y) \Rightarrow B(\varepsilon_y B)$ n'appartient pas à $\Delta$, alors $\Delta'$ est une preuve de $\varphi$ dans $CP_\varepsilon^*$ et l'on conclut par récurrence. Si au contraire $\exists y \, B(y) \Rightarrow B(\varepsilon_y B)$ appartient à $\Delta$, alors $\Delta'$ est une preuve de $\varphi$ dans $CP_\varepsilon^*$ à partir de $\exists y \, B \Rightarrow B(a)$. On a donc dans $CP_\varepsilon^*$ une preuve $\Delta''$ de la formule

(1) $\quad (\exists y \, B(y) \Rightarrow B(a)) \Rightarrow \varphi$

qui est sans symbole $\varepsilon$ par hypothèse. Par hypothèse de récurrence, (1) est dérivable dans CP. Mais $B(a) \Rightarrow \varphi$ est une conséquence de (1) dans CP. On a donc $\vdash_{CP} B(a) \Rightarrow \varphi$ et par conséquent $\vdash_{CP} \exists y \, B(y) \Rightarrow \varphi$. Mais $\varphi$ est une conséquence dans CP de (1) et de $\exists y \, B(y) \Rightarrow \varphi$. Donc $\vdash_{CP} \varphi$. On généralise facilement au cas $\Sigma \vdash_{CP_\varepsilon^*} \varphi$.

---

[25] Cf. par exemple Leisenring [1969].

## 2. Le rapport à la théorie des ensembles

Le fait que $CP_\varepsilon$ soit en logique classique une extension inessentielle de CP est dû au fait que l'on ne considère que des axiomes sans symbole $\varepsilon$. Si l'on se permet d'introduire des $\varepsilon$-termes *dans les axiomes* des théories, alors on peut obtenir des extensions *effectives*. Tel est le cas par exemple de la théorie des ensembles *où l'axiome du choix devient dérivable*. Cet axiome dit que si $X = \{X_\alpha\}_{\alpha \in I}$ est une famille d'ensemble non vides $X_\alpha$ indexés par un ensemble $I$, il existe *un ensemble Y* comprenant exactement *un* élément de $X_\alpha$ pour chaque $\alpha$.

Considérons alors la collection $C = \{\varepsilon_z(z \in X_\alpha)\}$ des $\varepsilon$-termes canoniquement associés aux $X_\alpha$. D'après *l'axiome de compréhension* disant que pour tout ensemble $Z$ et toute propriété $P$ définie sur $Z$, il existe un ensemble $V$ constitué des éléments de $Z$ satisfaisant $P$, la collection $C$ est un ensemble. Il suffit en effet de prendre pour $Z$ l'union

$$Z = \bigcup_{\alpha \in I} X_\alpha$$

(qui est un ensemble d'après l'axiome de réunion de la théorie des ensembles) et pour $P$ la propriété $P$ définie par

$$P(x) \equiv \exists \alpha \big(\alpha \in I \wedge x = \varepsilon_z(z \in X_\alpha)\big).$$

$C$ étant un ensemble, on peut le prendre pour ensemble $Y$.

L'opérateur $\varepsilon$ peut rendre beaucoup d'autres services en théorie des ensembles.[26] Par exemple, comme y a insisté Carnap, il permet très facilement de définir le concept de *cardinal* car celui-ci possède exactement son degré de généralité. Soit $x \sim y$ la relation d'équipotence entre ensembles et $\#x$ le cardinal de $x$. Par définition des cardinaux comme classes d'équivalence pour l'équipotence, on veut que la formule : $\forall x \forall y \, (x \sim y \Leftrightarrow \#x = \#y)$ soit valide. La définition la plus simple consiste à prendre : $\#x \equiv \varepsilon_y(y \sim x)$. Cette procédure est d'ailleurs valable pour toute relation d'équivalence. Elle permet de sélectionner systématiquement des représentants dans les classes d'équivalence, ce qui, ainsi que l'ont souligné Beth et Carnap est fondamental pour comprendre *les définitions par abstraction*.

Mais l'opérateur $\varepsilon$ peut rendre des services plus subtils, par exemple pour comprendre ce qu'est une formule *collectivisante*. Soit $B(x)$ une formule. On dit qu'elle est collectivisante si les éléments qui y satisfont forment un ensemble $X_B$, autrement dit si la formule $\exists X \forall x \, ((x \in X) \Leftrightarrow B(x))$ — qui est de la forme $\exists X \, F(X)$ — est dérivable. Mais l'on peut, indépendamment de toute démonstration d'existence, introduire *le concept* de $X_B$ en

---

[26] Cf. par exemple Leisenring [1969].

posant

$$X_B = \varepsilon_X F(X) = \varepsilon_X(\forall x\,((x \in X) \Leftrightarrow B(x))).$$

Dire que $B$ est collectivisante revient alors à affirmer *le schéma de Church* : $\forall x\,((x \in X_B) \Leftrightarrow B(x))$ qui n'est qu'un exemple d'ε-formula du $2^{\text{ème}}$ ordre $\exists X\,F(X) \Leftrightarrow F(\varepsilon_F) \Leftrightarrow F(X_B)$.

## V.     L'ε-CALCULUS ET LA CONSISTANCE DE L'ARITHMETIQUE

Toute la "machinerie" de l'ε-calcul a été mise en place par Hilbert pour pouvoir transformer des preuves arithmétiques quelconques (non finitistes) en preuves *finitistes* purement *combinatoires*. Son idée de base était d'utiliser le calcul $CP_\varepsilon$ à la place du calcul des prédicats classique CP puis d'éliminer ensuite des preuves, au moyen d'une *méthode de substitution numérique,* les "formules critiques" instantiant l'ε-formula (ε) $F(t) \Rightarrow F(\varepsilon_F)$.

Soit $P$ une preuve et $E$ l'ensemble (fini) des formules critiques (ε) intervenant dans $P$. On effectue des substitutions $S$ en remplaçant les ε-termes saturés par des *valeurs numériques*. Si toutes les formules critiques deviennent vraies sous $S$ on dit que $S$ résout $P$. L'idée est alors de trouver une $S$ résolvante par approximations successives en partant de $S_0$ qui remplace tous les ε-termes par 0 puis en corrigeant pas à pas les valeurs numériques attribuées. Tout le problème est de montrer qu'un tel processus se termine.

La principale difficulté est de savoir comment attribuer des valeurs numériques à des ε-termes *enchâssés*. Supposons que $P$ utilise une formule critique simple $B(y) \Rightarrow B(\varepsilon_B)$ et également une formule critique complexe $A(x,\varepsilon_B) \Rightarrow A(\varepsilon_A,\varepsilon_B)$ où $\varepsilon_A = \varepsilon_x A(x,\varepsilon_B)$. $\varepsilon_B$ doit être évalué avant $\varepsilon_A$. Mais si l'on doit utiliser aussi la formule critique $B(\varepsilon_A) \Rightarrow B(\varepsilon_B)$ alors on rencontre un problème. En effet supposons que $B(0)$ soit faux et commençons par évaluer $\varepsilon_B = 0$. Si l'on change ensuite l'évaluation de $\varepsilon_x A(x,0)$ de 0 à $n$, on aura, d'après la formule $B(\varepsilon_A) \Rightarrow B(\varepsilon_B)$, la formule $B(n) \Rightarrow B(0)$ qui est fausse si $B(n)$ est vraie (puisque $B(0)$ est faux par hypothèse). Il faudra donc évaluer $\varepsilon_B = n$. Mais alors c'est $\varepsilon_x A(x,n)$ qui ne sera plus forcément vraie.

Les ε-termes du deuxième ordre introduisent des difficultés encore plus insurmontables. Par exemple considérons l'ε-axiome du deuxième ordre

$$A_a(F(a)) \Rightarrow A_a(\varepsilon_F A_b(F(b))(a))$$

où la notation $A_a$ signifie que $a$ n'est pas une variable libre dans $A$ et n'intervient qu'à travers

25l'argument fonctionnel *F*. Dans la mesure où $\varepsilon_F A_b(F(b))$ est la fonction définie par *A*, cet axiome est l'analogue d'un axiome de compréhension. L'introduction des ε-termes du deuxième ordre conduit donc à restreindre les axiomes de compréhension dans le système d'Ackermann, ce qui rend ce dernier *prédicatif* comme l'a remarqué von Neumann en 1927 (cf. Zach, 2002).

Mais, quoi qu'il en soit, la consistance de l'arithmétique du premier ordre puis de l'arithmétique primitive récursive du second ordre ont été démontrées par Ackermann dans sa grande thèse de 1924 *Begründung des 'tertium non datur' mittels der Hilbertschen Theorie des Widerspruchsfreiheit*, puis clarifié ensuite en 1940 à la suite de la preuve de Gentzen de 1936, Ackermann montrant qu'une induction transfinie jusqu'à l'ordinal $\omega^{\omega^\omega}$ est nécessaire (cf. les travaux de Grigori Mints et le survey de Richard Zach (2002)).

Il s'agissait donc d'abord de démontrer de façon finitiste la consistance de l'arithmétique formelle, c'est-à-dire du système d'axiomes 𝕹 constitué :

(i)     du système $𝕹^0$ axiomatisant à la Peano la structure d'ordre et la structure algébrique (+, .) de ℕ et

(ii)    du *schéma d'induction* :

(I)     $A(0) \Rightarrow (\forall x\, (A(x) \Rightarrow A(x+1)) \Rightarrow \forall x\, A(x))$ ou si l'on préfère
         $(A(0) \wedge \forall x\, (A(x) \Rightarrow A(x+1))) \Rightarrow \forall x\, A(x)$

où *A*(*x*) est une formule quelconque. Il existe une preuve finitiste de la consistance de $𝕹^0$ et c'est donc le schéma (I) du raisonnement par récurrence qui fait problème.

Hilbert utilise un ε-calcul où l'opérateur ε s'interprète comme un opérateur de *minimum*. Si *F*(*x*) est un prédicat sur ℕ tel que $X_F \neq \emptyset$ alors $\varepsilon_F$ dénote *le plus petit* élément de $X_F$ (et si $X_F = \emptyset$ alors $\varepsilon_F = 0$). Cette interprétation est la seule qui satisfasse le schéma d'axiomes :

$F(a) \Rightarrow \varepsilon_F \leq a.$

Pour démontrer la consistance de 𝕹 il faut et il suffit de démontrer que tout énoncé *élémentaire* (i.e. sans quantificateurs et sans symbole ε) qui est un $CP_\varepsilon$ théorème de 𝕹 devient numériquement vrai une fois interprété dans ℕ. Le problème posé par le schéma d'axiomes (I) est que les axiomes n'y sont pas sous forme prénexe et que leur résolution de Skolem conduit à ajouter au système 𝕹´ que l'on voudrait substituer à 𝕹 toutes les fonctions de Skolem, ce qui complique énormément les choses. Dans leurs *Grundlagen*, Hilbert et Bernays définissent alors le système d'axiomes $𝕹_\varepsilon$ obtenu à partir des axiomes de 𝕹⁰ (qui sont tous du type $\forall x\, A(x)$ ou $\forall x \forall y\, B(x,y)$) en prenant leurs matrices et en y substituant les variables libres $x$ et $y$ par des symboles de constantes $s$ et $t$ et en ajoutant les schémas d'axiomes :

(E₁)   $\neg A(\varepsilon_A) \Rightarrow \neg A(t)$

(E₂)   $\varepsilon_A = t+1 \Rightarrow \neg A(t)$.

Ils démontrent alors que l'on peut substituer $𝕹_\varepsilon$ à $𝕹$ et que l'induction est démontrable dans $𝕹_\varepsilon$.

*Théorème* (preuve de l'induction). Soit $B$ une formule inductive :
$$B \equiv A(0) \Rightarrow (\forall x\, (A(x) \Rightarrow A(x+1))) \Rightarrow \forall x\, A(x)).$$
Alors $𝕹_\varepsilon \vdash_{CP_\varepsilon} B$. Qui plus est, si $B$ est *propre*, autrement dit si $x$ ne possède pas dans $A(x)$ d'occurrence libre qui soit dans la portée d'un $\varepsilon$-symbole, alors $𝕹_\varepsilon \vdash_{CP^*_\varepsilon} B$.

*Preuve.*[27] Notons $s$ l'$\varepsilon$-terme $\varepsilon_{\neg A}$. Soit $\Sigma$ le système d'axiomes : $\Sigma = 𝕹_\varepsilon \cup \{A(0), \forall x\, (A(x) \Rightarrow A(x+1))\}$. Comme $\forall x\, A(x) \Leftrightarrow A(\varepsilon_{\neg A})$, on veut en fait montrer que l'on a .

Posons $t = g(s) = \begin{cases} 0 & \text{si } s = 0 \\ s-1 & \text{si } s \geq 1 \end{cases}$.

On a la démonstration dans $CP_\varepsilon$ :

(1)   $\Sigma \vdash (s = 0) \vee (s = t+1)$ (par définition de $t$);

(2)   $\Sigma \vdash A(0)$ (par définition de $\Sigma$);

(3)   $\Sigma \vdash ((s = 0) \wedge A(0)) \Rightarrow A(s)$ (trivial);

(4)   $\Sigma \vdash (s = 0) \Rightarrow A(s)$ (d'après (3) et (2));

Donc si $s = 0$, on a bien $\Sigma \,O\, A(s)$.

Si en revanche $s \geq 1$, on a :

(5)   $\Sigma \vdash \forall x\, (A(x) \Rightarrow A(x+1))$ (par définition de $\Sigma$);

(6)   $\Sigma \vdash A(t) \Rightarrow A(t+1)$ (d'après (5));

(7)   $\Sigma \vdash (s = t+1) \wedge A(t+1) \Rightarrow A(s)$ (trivial);

Mais $s = \varepsilon_{\neg A}$ par définition et $s = t+1$ par hypothèse (car $s \geq 1$). Donc d'après l'axiome (E₂) appliqué à $\neg A$, on a :

---

[27] Nous donnons ici cette preuve (élémentaire) pour que le lecteur non mathématicien puisse avoir une toute petite idée des techniques hilbertiennes.

(8)     $\Sigma \vdash (s = t+1) \Rightarrow \neg\neg A(t)$;

donc par (8), $\neg\neg A(t) = A(t)$, (6) et (7) :

(9)     $\Sigma \vdash (s = t+1) \Rightarrow A(s)$;

donc par (1), (4) et (9) :

(10)    $\Sigma \vdash A(s)$.

Ce théorème montre que $\mathfrak{N}_\varepsilon$ est une extension inessentielle de $\mathfrak{N}$. On utilise ensuite les ε-théorèmes pour montrer que si *B* est un théorème élémentaire de $\mathfrak{N}$ alors il en existe une preuve dans CE. Il ne reste donc plus qu'à démontrer la

*Conjecture*. Si *B* est une formule élémentaire démontrable dans CE, alors *B* est numériquement vraie.

C'est cette conjecture que Hilbert, Bernays et Ackermann n'ont pu démontrer et cela parce qu'en fait elle est invalidée par le théorème d'incomplétude de Gödel.

Toutefois, Hilbert et Bernays ont abouti à un résultat partiel. On appelle *arithmétique restreinte*, la théorie $\mathfrak{N}^r$ où l'application du schéma d'induction (I) est restreinte aux formules *propres*. L'élimination de l'impropreté permet de démontrer la conjecture et de fournir une preuve finitiste de la consistance de $\mathfrak{N}^r$.

## VI.    LA NATURE INTENSIONNELLE DES ε-TERMES

Nous avons déjà signalé plusieurs fois la nature *intensionnelle* des ε-termes et leur double interprétation *de dicto* et *de re*. Elle provient de leur structure de symbole-index. Rappelons que les symboles-index (au sens de Peirce, ou encore les "expressions indexicales" au sens de Bar-Hillel) correspondent dans la classification de Jakobson aux unités linguistiques pour lesquelles et par lesquelles le code renvoie au message. Se situant à la charnière du sens et de la dénotation (du *Sinn* et de la *Bedeutung*), ils sont à l'origine des *contextes obliques* (aussi dits "opaques"). Les exemples linguistiques standard en sont les démonstratifs (celui-ci, etc.), les déictiques (ici, maintenant, etc.), et les shifters (je, tu, etc.). Ils sont de nature pragmatique, font intervenir l'instance de l'énonciation et leur dénotation n'est définissable que relativement au contexte. Qui plus est, comme le notait Emile Benvéniste, ce sont des formes linguistiques qui, contrairement aux termes nominaux, ne réfèrent qu'à des individus. Ils ne sont pas abstraits mais concrets et intensionnellement décomplétés.

En général cette essence pragmatique n'est pas véritablement analysée comme telle, même chez des linguistes comme Morris car elle reste réduite à la dimension empirique du

28rapport qu'un énoncé entretient avec la situation énonciative. Ce n'est que relativement récemment que certains linguistes, par exemple Jean-Pierre Desclés, en sont venus à la considérer comme linguistiquement primitive :

> "On conçoit mal un modèle général traitant des langues naturelles où, ce qui apparaît comme essentiellement langagier et général, serait traité dans la composante pragmatique ajustant les composantes syntaxiques et sémantiques aux données empiriques et aux concepts "intuitifs", la pragmatique étant ajoutée pour "rendre compte" des quelques phénomènes déclarés être auxiliaires, alors qu'ils semblent présents dans toutes les langues et sont, selon nous, beaucoup plus primitifs que la plupart des autres phénomènes et peut-être même au sens de l'activité de langage."[28]

La logique hilbertienne est une logique pragmatique de la généricité qui conduit à associer à tout concept une entité intensionnelle idéale structurée comme un symbole-index. Ce lien entre indexicalité et typicalité lui confère une grande portée sémio-linguistique et philosophique. Une interprétation adéquate des ε-termes doit la faire intervenir d'une façon ou d'une autre, par exemple en faisant que la dénotation $a$ de $\varepsilon_F$ soit relative à un *contexte* de sélection. Von Heusinger et Wilfrid Meyer Viol ont beaucoup insisté sur ce point ces dernières années.

Mais il y a déjà longtemps que Melvin Fitting [29] a montré comment on pouvait développer une interprétation intensionnelle des $\varepsilon$-termes dans le cadre d'une logique kripkéenne des mondes possibles. On introduit un ensemble $\mathfrak{M}$ de "mondes possibles $M_i$ où un monde se trouve distingué comme "monde réel" $M_0$. On se donne en plus sur $\mathfrak{M}$ une relation $\mathfrak{R}$ d'*accessibilité* (réflexive) $M_1 \mathfrak{R} M_2$. On se donne enfin des assignations $V_{M_i}(p)$ assignant à toute proposition $p$ une valeur de vérité dans $M_i$. Si $V_{M_i}(p)$ est vraie, on note $\vDash_{M_i} p$.

Lorsque l'on passe au calcul des prédicats CP, on peut alors naturellement étendre ce genre de sémantique modale. On se donne d'abord des assignations cohérentes de valeurs de vérité pour les propositions, la cohérence signifiant (pour $\mathfrak{R}$ transitive par exemple) que :

(i)     $\vDash_M p \wedge q$ ssi $\vDash_M p$ et $\vDash_M q$,

---

[28] Tel est le cas par exemple de Jean-Pierre Desclés [1990]. Il a consacré de nombreuses études à la formalisation des structures de l'énonciation. Ses recherches sur la cognition l'ont d'ailleurs conduit à remettre en chantier le formalisme hilbertien dans un cadre de pensée original.

[29] Fitting [1970].

(ii)     ⊨ $_M$ $p \Rightarrow q$ ssi pour tout $M'$ tel que $M \mathfrak{R} M'$ on a non ⊨ $_{M'} p$ ou ⊨ $_{M'} q$,

(iii)    ⊨ $_M$ $\neg p$ ssi pour tout $M'$ tel que $M \mathfrak{R} M'$ on a non ⊨ $_{M'} p$, etc.

On dit alors qu'une proposition complexe $\varphi$ est *valide* si ⊨ $_{M_0} \varphi$ pour *toute* assignation cohérente.

Pour pouvoir passer à des formules contenant des constantes et des variables, on suppose que chaque "monde" $M$ possède un domaine de constantes $D(M) \neq \emptyset$ tel que si $M \mathfrak{R} M'$ alors $D(M) \subseteq D(M')$. Si $F(x_1, ..., x_n)$ est un prédicat $n$-aire on lui associe dans chaque monde $M$ un sous-ensemble $V_M(F)$ de $D(M)^n$ avec la condition que si $M \mathfrak{R} M'$ alors $V_M(F) \subseteq V_{M'}(F)$. Soit $\mathscr{U}$ l'univers des constantes (l'union des $D(M)$) et soit $\varphi$ une formule contenant les variables $x_1, ..., x_n$. Considérons une assignation dans $\mathscr{U}$, $a_1, ..., a_n$, des valeurs de ces variables. On pose ⊨ $_M \varphi$ ssi les constantes de $\varphi$ appartiennent à $D(M)$ et $(a_1, ..., a_n) \in V_M(\varphi)$. D'où les définitions évidentes :

(i)     ⊨ $_M \exists y\, \varphi(x,y)$ pour une assignation $x \rightarrow a$ ssi il existe $b \in D(M)$ tel que ⊨ $_M \varphi(a,b)$,

(ii)    ⊨ $_M \forall y\, \varphi(x,y)$ pour une assignation $x \rightarrow a$ ssi pour tout $M'$ tel que $M \mathfrak{R} M'$ et pour toute assignation $y \rightarrow b \in D(M')$, on a ⊨ $_{M'} \varphi(a,b)$.

Dans un cadre de ce genre l'on peut alors, de façon naturelle, développer une interprétation intensionnelle de l'opérateur $\varepsilon$, par exemple pour le système modal S$_4$.[30] Les difficultés liées aux interprétations non intensionnelles des $\varepsilon$-termes deviennent criantes dans une sémantique modale à la Kripke. En effet, si $F(x)$ est un prédicat (unaire), $\varepsilon_F$ est un symbole d'individu. Mais dans un modèle kripkéen de S$_4$ l'existentielle $\exists x\, F(x)$ peut être valide dans deux mondes possibles sans qu'il n'existe pourtant aucune constante $c$ satisfaisant $F(c)$ dans les deux mondes. Il faut donc traiter $\varepsilon_F$ comme *une fonction sur les mondes possibles* telle que si ⊨ $_M \exists x\, F(x)$ alors la valeur de $\varepsilon_F$ dans $M$ est une constante $c$ telle que ⊨ $_M F(c)$. Les $\varepsilon$-termes ne sont donc ni des variables, ni des constantes mais des entités *relatives aux mondes* ("world-dependent") *dont la dénotation varie avec les mondes*. Leur caractère "mixte" qui avait été admirablement dégagé par Lautman (cf. plus haut) est bien un caractère intensionnel.

Fitting montre en fait qu'il existe une fonction $f_F$ telle que :

(i)     le domaine $dom(f_F)$ de $f_F$ est l'ensemble des mondes $M \in \mathfrak{M}$ sur lesquels $F$ est définie (i.e. tels que les constantes de $F$ appartiennent à $D(M)$);

(ii)    si $M \in dom(f_F)$ alors $f_F(M) \in M$;

(iii)   si ⊨ $_M \exists x\, F(x)$ alors ⊨ $_M F(x \mid f_F(M))$.

Autrement dit, $f_F$ est une fonction de choix qui, pour chaque monde $M$ sur lequel $F$ est

---

[30] Cf. Fitting [1970].

définie et où ∃*x F(x)* est valide, sélectionne le référent de $\varepsilon_F$.

Ce qui est particulièrement intéressant dans cette construction est qu'elle mime logiquement de nombreux aspects sémiotiques classiquement opposés à la sémantique vériconditionnelle classique.

(i)     $\varepsilon_F$ est un symbole d'individu mais possède une structure interne et une signification dans la mesure où il est canoniquement associé à *F*.

(ii)    Il ne dénote pas à travers un renvoi symbolique. Il sélectionne en fonction de sa structure interne des éléments dans des mondes possibles. Cette sélection pragmatique dépend d'une instance de choix intensionnelle.

(iii)   Le référent sélectionné est une sorte d'objet dynamique à la Peirce qui est indexicalement sélectionné en tant que représenté d'une certaine façon, i.e. en tant que représenté par *F*.

(iv)    L'existence devient par conséquent indépendante d'un rapport dénotationnel au monde et s'identifie à une consistance de la signification.

(v)     Il peut du coup y avoir une inférence concernant l'existence. L'existence se disjoint de l'ontologie et acquiert le statut d'un prédicat.

## VII.   L'ε-CALCULUS ET LA LOGIQUE INTUITIONNISTE

Des travaux plus récents ont étudié l'opérateur ε dans le cadre de la logique intuitionniste où l'implication $F \Rightarrow \neg\neg F$ est valide (car si *a* satisfait *F*, alors *a* est un contre exemple explicite à ¬*F(x)* et on a donc ¬¬*F(x)*) sans que le soit pour autant l'implication réciproque :

(¬¬)    $\neg\neg F \Rightarrow F$.

La différence entre *F* et ¬¬*F* se voit très bien dans la sémantique intuitionniste standard où les valeurs de vérité constituent une algèbre de Heyting et non plus une algèbre de Boole. On y considère que les extensions des prédicats, au lieu d'être des ensembles quelconques, sont des ouverts d'espaces topologiques. Si $X = X_F$ est l'extension (ouverte) de *F*, alors celle de ¬*F* n'est pas le complémentaire $\complement X$ de *X* (qui est fermé) mais son *intérieur* $\overset{\circ}{\complement}X$. Cela implique le viol du tiers exclu *F*∨¬*F* car pour les éléments du bord ∂*X* de *X* on n'a *ni F(x) ni* ¬*F(x)*. Comme le complémentaire de $\overset{\circ}{\complement}X$ est la fermeture $\overline{X}$ de *X*, on voit que ¬¬*F* correspond à $\overset{\circ}{\overline{X}}$, l'intérieur de la fermeture de *X*. On a évidemment (car *X* est ouvert) $X \subset \overset{\circ}{\overline{X}}$. Les prédicats pour lesquel $F = \neg\neg F$ sont dits *décidables*.

La logique intuitionniste permet de résoudre de façon particulièrement efficace un certain nombre de paradoxes, par exemple ceux des infinitésimales leibniziennes évoqué plus haut. En effet, on peut montrer que *dx* devient un symbole dénotant dans la théorie de ℝ le sous-ensemble $\neg\neg\{0\}$ des réels non non nuls.[31] En logique classique $\neg\neg\{0\}=\{0\}$ et il n'existe donc pas d'infinitésimales. Mais en logique intuitionniste $\neg\neg\{0\}\neq\{0\}$ et il existe donc des infinitésimales.

On peut développer une version intuitionniste $CPI_\varepsilon$ de l'$\varepsilon$-calculus hilbertien. On garde l'axiome d'équivalence définissant $\varepsilon$ :

($\exists$)     $\exists x\, F(x) \Leftrightarrow F(\varepsilon_F)$.

Mais il n'implique plus pour les universelles l'équivalence:

($\forall$)     $\forall x\, F(x) \Leftrightarrow F(\varepsilon_{\neg F})$.

On a évidemment $\forall x\, F(x) \Rightarrow F(\varepsilon_{\neg F})$ et donc $\neg F(\varepsilon_{\neg F}) \Rightarrow \neg \forall x\, F(x)$, i.e., d'après ($\exists$), $\exists x\, \neg F(x) \Rightarrow \neg \forall x\, F(x)$. Mais l'on n'a pas l'implication réciproque $\neg \forall x\, F(x) \Rightarrow \neg F(\varepsilon_{\neg F})$ ou $\neg \forall x\, F(x) \Rightarrow \exists x\, \neg F(x)$ (ce qui est un aspect de la non validité du tiers exclu). Comme nous l'avons vu plus haut on n'a que le principe de Markov affirmant le tiers exclu pour les prédicats décidables.

Contrairement à ce qui se passe dans le cas classique, $CPI_\varepsilon$ *n'est pas* une extension conservative de CPI. Il est strictement plus riche. Par exemple on peut y construire effectivement des fonctions *discontinues* sur ℝ alors que dans la théorie intuitionniste de ℝ toutes les fonctions sont continues.

Dans une fort belle étude sur l'$\varepsilon$-calculus intuitionniste, le logicien John Bell (qui a par ailleurs donné une version de l'ε-calculus dans le cadre de la théorie des topoï) a montré que l'écart entre logiques respectivement classique et intuitionniste était essentiellement lié au caractère *intensionnel* des $\varepsilon$-termes.

> "Since the $\varepsilon_A$ have been introduced *intensionally*, that is by the form of the defining predicate *A*, we do not yet possess a useful sufficient *identity* condition." [32]

On peut rétablir l'extensionalité avec l'axiome d'identité d'Ackermann ($\varepsilon_2$)

($\varepsilon_2$)     $\forall x\, (F(x) \Leftrightarrow G(x)) \Rightarrow \varepsilon_F = \varepsilon_G$

car alors

---

[31] Pour une introduction à ces questions, et en particulier aux résultats de Penon, cf. Petitot [1999].

[32] Bell [1993], p. 6.

"the identity of $\varepsilon_A$ is completely determined by the extension of A."

John Bell démontre alors le théorème que *l'axiome d'Ackermann implique le tiers exclu $F \vee \neg F$. Ce résultat philosophiquement profond montre que la logique classique est une conséquence du fait que les prédicats peuvent être exemplifiés de façon extensionnelle.*

"In short, not existence, but extensional existence, yields classical logic." [33]

John Bell montre ensuite que l'$\varepsilon$-calculus intuitionniste $CPI_\varepsilon$ est "sain" (ou "adéquat"), i.e. $\Sigma \vdash_{CP_\varepsilon} \varphi$ implique $\Sigma \vDash_{CP_\varepsilon} \varphi$. Mais, contrairement au $CP_\varepsilon$ classique, *il n'est pas complet* car, par exemple, la formule $\varepsilon_x(x = x) = \varepsilon_x(x \neq x)$ y est valide sans y être pour autant dérivable.

Enfin John Bell montre qu'il est *impossible* de dériver dans $CPI_\varepsilon$ l'inférence $\neg \forall x\, F(x) \Rightarrow \exists x\, \neg F(x)$ ainsi que le tiers exclu $F \vee \neg F$. Pour démontrer ces assertions, le principe d'extensionalité d'Ackermann se révèle absolument nécessaire.

**BIBLIOGRAPHIE**

---

[33] *Ibid.*, p. 8.